\documentclass[12pt,leqno]{article}

\newtheorem{theorem}{Theorem}
\newtheorem{lemma}[theorem]{Lemma}
\newtheorem{proposition}[theorem]{Proposition}
\newtheorem{definition}[theorem]{Definition}
\newtheorem{corollary}[theorem]{Corollary}

\newcommand{\begintheorem}{\addtocounter{equation}{1}\begin{theorem}}
\newcommand{\beginlemma}{\addtocounter{equation}{1}\begin{lemma}}
\newcommand{\beginproposition}{\addtocounter{equation}{1}\begin{proposition}}
\newcommand{\begindefinition}{\addtocounter{equation}{1}\begin{definition}}
\newcommand{\begincorollary}{\addtocounter{equation}{1}\begin{corollary}}

\begin{document}

\title{Potpourri, 4}

\author{Stephen William Semmes	\\
	Rice University		\\
	Houston, Texas}

\date{}

\maketitle

\tableofcontents

\section{Random polynomials}
\label{random polynomials}
\setcounter{equation}{0}

	Fix a positive integer $n$, and let us write
$\mathcal{B}_{n+1}$ for the set of binary strings of length $n + 1$
written multiplicatively.  To be more precise, an element $s$ of
$\mathcal{B}_{n+1}$ is an $(n+1)$-tuple $s = (s_0, \ldots, s_n)$ such
that each $s_j$ is either $1$ or $-1$.  Put $r_j(s) = s_j$ for each
$j$, the $j$th Rademacher function on $\mathcal{B}_{n+1}$.  Let $a_0,
\ldots, a_n$ be complex numbers, so that we get a polynomial
\begin{equation}
	p(z) = \sum_{j=0}^n a_j \, z^j,
\end{equation}
where $z^j$ is interpreted as being equal to $1$ when $j = 0$, as
usual.  For each $s \in \mathcal{B}_{n+1}$ consider the polynomial
\begin{equation}
	p_s(z) = \sum_{j=0}^n a_j \, r_j(t) \, z^j.
\end{equation}

	Let us write ${\bf T}$ for the unit circle in the complex
plane, which is to say the set of complex numbers $z$ with $|z| = 1$.
For each positive integer $j$ we have that
\begin{equation}
	\int_{\bf T} z^j \, |dz| = 0,
\end{equation}
and hence also
\begin{equation}
	\int_{\bf T} \overline{z}^j \, |dz| = 0,
\end{equation}
where $\overline{z}$ denotes the complex conjugate of $z$.  Of course
\begin{equation}
	\frac{1}{2 \pi} \int_{\bf T} |dz| = 1.
\end{equation}
As usual, this leads to
\begin{equation}
	\frac{1}{2 \pi} \int_{\bf T} |p(z)|^2 \, |dz| = \sum_{j=0}^n |a_j|^2.
\end{equation}
Similarly,
\begin{equation}
	\frac{1}{2 \pi} \int_{\bf T} |p_s(z)|^2
		= \sum_{j=0}^n |a_j|^2
\end{equation}
for all $s \in \mathcal{B}_{n+1}$.

	Now consider
\begin{equation}
	2^{-n-1} \sum_{s \in \mathcal{B}_{n+1}} 
		\frac{1}{2 \pi} \int_{\bf T} |p_s(z)|^4 \, |dz|.
\end{equation}
By interchanging the order of summation and integration we can rewrite
this as
\begin{equation}
	\frac{1}{2 \pi} \int_{\bf T} 2^{-n-1} 
		\sum_{s \in \mathcal{B}_{n+1}} |p_s(z)|^4 \, |dz|.
\end{equation}
There is a positive real number $C$ such that
\begin{equation}
	2^{-n-1} \sum_{s \in \mathcal{B}_{n+1}} 
			\biggl|\sum_{j=0}^n b_j \, r_j(s) \biggr|^4
		\le C \, \bigg(\sum_{j=0}^n |b_j|^2 \bigg)^2
\end{equation}
for any complex numbers $b_0, \ldots, b_n$, as one can show by
expanding
\begin{equation}
	\biggl|\sum_{j=0}^n b_j \, r_j(s)\biggl|^4
\end{equation}
as a quadruple sum and then summing over $s$ first.  It follows that
\begin{equation}
	2^{-n-1} \sum_{s \in \mathcal{B}_{n+1}} 
		\frac{1}{2 \pi} \int_{\bf T} |p_s(z)|^4 \, |dz|
	\le C \, \bigg(\sum_{j=0}^n |a_j|^2 \bigg)^2.
\end{equation}

	More generally, suppose that $m$ is a positive integer,
and consider
\begin{equation}
	2^{-n-1} \sum_{s \in \mathcal{B}_{n+1}} 
		\frac{1}{2 \pi} \int_{\bf T} |p_s(z)|^{2m} \, |dz|.
\end{equation}
As before there is a positive real number $C(m)$ such that
\begin{equation}
	2^{-n-1} \sum_{s \in \mathcal{B}_{n+1}} 
			\biggl|\sum_{j=0}^n b_j \, r_j(s) \biggr|^{2m}
		\le C(m) \, \bigg(\sum_{j=0}^n |b_j|^2 \bigg)^m
\end{equation}
for any complex numbers $b_0, \ldots, b_n$.  By interchanging the sum
with the integral and then applying this fact one obtains that
\begin{equation}
	2^{-n-1} \sum_{s \in \mathcal{B}_{n+1}} 
		\frac{1}{2 \pi} \int_{\bf T} |p_s(z)|^4 \, |dz|
		\le C(m) \, \bigg(\sum_{j=0}^n |a_j|^2 \bigg)^m.
\end{equation}

\section{Vector-valued functions}
\label{vector-valued functions}
\setcounter{equation}{0}

	Let $E$ be a nonempty finite set, and let $V$ be a real or
complex vector space.  The vector space of functions on $E$ with
values in $V$ will be denoted $\mathcal{F}(E, V)$.

	In particular, $\mathcal{F}(E, {\bf R})$, $\mathcal{F}(E, {\bf
C})$ denote the spaces of real and complex-valued functions on $E$,
respectively.  These are commutative algebras with respect to
pointwise multiplication of functions.

	Actually, we can think of $\mathcal{F}(E, V)$ as a module over
$\mathcal{F}(E, {\bf R})$ or $\mathcal{F}(E, {\bf C})$, according to
whether $V$ is a real or complex vector space.  In other words, a
function on $E$ with values in $V$ can be multiplied pointwise by a
scalar-valued function in a way that is compatible with addition and
scalar multiplication.  For that matter, scalar multiplication amounts
to the same thing as multiplication by constant scalar-valued
functions.

	We can also think of $\mathcal{F}(E, V)$ as a module over the
algebra of linear transformations on $V$, where a linear
transformation on $V$ induces a linear transformation on $V$-valued
functions on $E$ by acting on the values pointwise.  The actions on
$\mathcal{F}(E, V)$ by pointwise multiplication by scalar-valued
functions on $E$, and by pointwise action by linear transformations on
$V$, obviously commute with each other.

	Suppose that $V$ is equipped with a norm $\|v\|_V$.
Thus $\|v\|_V$ is a nonnegative real number for each $v \in V$
which is equal to $0$ if and only if $v = 0$,
\begin{equation}
	\|\alpha \, v \|_V = |\alpha| \, \|v\|_V
\end{equation}
for all real or complex numbers $\alpha$, as appropriate, and all $v
\in V$, and
\begin{equation}
	\|v + w\|_V \le \|v\|_V + \|w\|_V
\end{equation}
for all $v, w \in V$.  This choice of norm on $V$ leads to a metric
$\|v - w\|_V$ on $V$, as usual.

	Let $p$ be given, $1 \le p \le \infty$.  If $f(x)$ is a real
or complex-valued function on $E$, put
\begin{equation}
	\|f\|_p = \bigg( \sum_{x \in E} |f(x)|^p \bigg)^{1/p}
\end{equation}
when $p < \infty$ and
\begin{equation}
	\|f\|_\infty = \max \{|f(x)| : x \in E\}.
\end{equation}
As is well-known, these define norms on the vector spaces of real and
complex-valued functions on $E$.  Moreover, if $1 \le p \le q \le
\infty$ and $f$ is a real or complex-valued function on $E$, then
\begin{equation}
	\|f\|_q \le \|f\|_p \le |E|^{(1/p) - (1/q)} \, \|f\|_q,
\end{equation}
where $|E|$ denotes the number of elements in $E$.

	If $f$ is a $V$-valued function on $E$, put
\begin{equation}
	\|f\|_{p, V} = \bigg( \sum_{x \in E} \|f(x)\|_V^p \bigg)^{1/p}
\end{equation}
when $1 \le p < \infty$ and
\begin{equation}
	\|f\|_{\infty, V} = \max \{\|f(x)\|_V : x \in E\}.
\end{equation}
These define norms on $\mathcal{F}(E, V)$.  Once again we have that
\begin{equation}
	\|f\|_{q, V} \le \|f\|_{p, E} \le |E|^{(1/p) - (1/q)} \, \|f\|_{q, V}
\end{equation}
when $1 \le p \le q \le \infty$ and $f \in \mathcal{F}(E, V)$.

	Let $V$ be a finite-dimensional real or complex vector space.
By a \emph{linear functional} on $V$ we mean a linear mapping from $V$
into the real or complex numbers, as appropriate.  One can add linear
functionals on $V$ and multiply them by scalars in the usual manner,
so that the space $V^*$ of linear functionals on $V$ is also a real or
complex vector space, as appropriate, called the dual of $V$.

	If $v_1, \ldots, v_n$ is a basis for $V$, so that every
element of $V$ can be expressed in a unique manner as a linear
combination of the $v_j$'s, then a linear functional $\lambda$ on $V$
is determined uniquely by the $n$ scalars $\lambda(v_1), \ldots,
\lambda(v_n)$, and for each choice of $n$ scalars there is a linear
functional on $V$ whose values on the basis vectors are those scalars.
In particular, $V^*$ is also finite-dimensional and has the same
dimension as $V$.

	Now suppose that $V$ is equipped with a norm $\|v\|_V$.  If
$\lambda$ is a linear functional on $V$, then there is a nonnegative
real number $k$ such that
\begin{equation}
	|\lambda(v)| \le k \, \|v\|_V
\end{equation}
for all $v \in V$.  Indeed, $V$ is isomorphic to ${\bf R}^n$ or ${\bf
C}^n$ as a vector space, where $n$ is the dimension of $V$, and it is
well-known that any norm on ${\bf R}^n$, ${\bf C}^n$ determines the
same topology as the standard Euclidean norm.

	As a result,
\begin{equation}
	\|\lambda\|_{V^*} = \sup \{|\lambda(v)| : v \in V, \ \|v\|_V \le 1\}
\end{equation}
is finite, and it is the same as the smallest nonnegative real number
which can be used as $k$ in the previous inequality.  This defines a
norm on the dual space $V^*$, which is the dual norm associated to
$\|\cdot \|_V$.

	Note that
\begin{equation}
	\|v\|_V = 
	   \sup \{|\lambda(v)| : \lambda \in V^*, \ \|\lambda\|_{V^*} \le 1\}.
\end{equation}
More precisely, $\|v\|_V$ is greater than or equal to $|\lambda(v)|$
for all linear functionals $\lambda$ on $V$ with dual norm less than
or equal to $1$ by definition, and there is such a linear functional
with $\lambda(v)$ equal to the norm of $v$ by famous duality results.

	As a basic example, let $E$ be a nonempty finite set, and let
us consider the vector space of real or complex-valued functions on
$E$.  If $h$ is a real or complex-valued function on $E$, then we get
a linear functional $\lambda_h$ on the vector space of real or
complex-valued functions on $E$, as appropriate, by setting
\begin{equation}
	\lambda_h(f) = \sum_{x \in E} f(x) \, h(x)
\end{equation}
for $f$ in the vector space.  Every linear functional on the vector
space of functions on $E$ arises in this manner.

	If $1 \le p, q \le \infty$ are conjugate exponents, in the
sense that $(1/p) + (1/q) = 1$, then H\"older's inequality implies
that
\begin{equation}
	|\lambda_h(f)| \le \|f\|_p \, \|h\|_q
\end{equation}
for all functions $f$, $h$ on $E$, where $\|f\|_p$, $\|h\|_q$ are as
in the previous section.  With respect to the norm $\|f\|_p$ on
functions on $E$, the dual norm of the linear functional $\lambda_h$
is therefore less than or equal to $\|h\|_q$, and in fact one can
check that it is equal to $\|h\|_q$.

	More generally we can consider $V$-valued functions on $E$.
If $h$ is a function on $E$ with values in $V^*$, then we can define a
linear functional $\lambda_h$ on the vector space of $V$-valued
functions on $E$ by saying that $\lambda_h(f)$ is obtained by applying
$h(x)$ as a linear functional on $V$ to $f(x)$ as an element of $V$
for each $x \in E$, and then summing over $x \in E$.  One can check
that every linear functional on $\mathcal{F}(E, V)$ occurs in this
way, so that the dual of $\mathcal{F}(E, V)$ can be identified with
$\mathcal{F}(E, V^*)$.

	Using H\"older's inequality it is easy to check that
$|\lambda_h(f)|$ is less than or equal to the product of $\|f\|_{p,
V}$ and $\|h\|_{q, V^*}$ for all $f \in \mathcal{F}(E, V)$ and $h \in
\mathcal{F}(E, V^*)$ when $1 \le p, q \le \infty$ are conjugate
exponents.  Furthermore, the dual norm of $\lambda_h$ with respect to
the norm $\|f\|_{p, V}$ on $\mathcal{F}(E, V)$ is exactly equal to
$\|h\|_{q, V^*}$.

	Let $V$ be a finite-dimensional real or complex vector space,
and let $E$ be a nonempty finite set.  If $f$ is a $V$-valued function
on $E$, then put
\begin{equation}
	\|f\|_{p, \nu} 
   = \sup \bigg\{\bigg(\sum_{x \in E} |\lambda(f(x))|^p \bigg)^{1/p} :
			\lambda \in V^*, \ \|\lambda\|_{V^*} \le 1 \bigg\}
\end{equation}
when $1 \le p < \infty$ and
\begin{equation}
	\|f\|_{\infty, \nu}
   = \sup \big\{ \max \{|\lambda(f(x))| : x \in E\} :
			\lambda \in V^*, \ \|\lambda\|_{V^*} \le 1 \big\}.
\end{equation}
One can check that these define norms on $\mathcal{F}(E, V)$.

	For each $v \in V$ we have that $\|v\|_V$ is equal to the
maximum of $|\lambda(v)|$ over $\lambda \in V^*$ with
$\|\lambda\|_{V^*} \le 1$, as discussed in the previous section.  It
is easy to see that
\begin{equation}
	\|f\|_{p, \nu} \le \|f\|_{p, V}
\end{equation}
when $1 \le p < \infty$, and that
\begin{equation}
	\|f\|_{\infty, \nu} = \|f\|_{\infty, V}.
\end{equation}
All of these norms reduce to $\|\cdot \|_V$ when $E$ has just one
element.

	We can also express $\|f\|_{p, \nu}$ as
\begin{equation}
	\quad \|f\|_{p, \nu} 
   = \sup \bigg\{\biggl|\sum_{x \in E} h(x) \, \lambda(f(x)) \biggr| :
	   \|h\|_q \le 1, \ \lambda \in V^*, \ \|\lambda\|_{V^*} \le 1 \bigg\}.
\end{equation}
Here $q$ is the conjugate exponent associated to $p$, so that $1 \le q
\le \infty$ and $(1/p) + (1/q) = 1$, and $h$ is a real or
complex-valued function on $E$, as appropriate.  Of course $\sum_{x
\in E} h(x) \, \lambda(f(x))$ is the same as $\lambda$ applied to
$\sum_{x \in E} h(x) \, f(x)$, and it follows that $\|f\|_{p, \nu}$ is
equal to the supremum of the $V$-norm of $\sum_{x \in E} h(x) \, f(x)$
over all scalar-valued functions $h$ on $E$ with $\|h\|_q \le 1$.

\section{Polynomials, continued}
\label{polynomials, continued}
\setcounter{equation}{0}

	Of course a polynomial on the complex plane has the form
\begin{equation}
\label{p(z) = a_n z^n + cdots + a_1 z + a_0}
	p(z) = a_n \, z^n + \cdots + a_1 \, z + a_0,
\end{equation}
where $a_0, \ldots, a_n$ are complex numbers.

	Alternatively one might consider polynomials in the real and
imaginary parts of $z$, which is equivalent to polynomials in $z$ and
$\overline{z}$.  A special case of this is given by linear
combinations of powers of $z$ and of powers of $\overline{z}$,
including constant terms, without products of positive powers of $z$
and $\overline{z}$.  One might also consider linear combinations of
powers of $z$ and of $z^{-1}$, including constant terms.  These
classes all define the same functions on the unit circle, where $|z|^2
= z \, \overline{z} = 1$.  Let us note that if $f(z)$ is one of these
more general kinds of polynomials, then there is a complex polynomial
$p(z)$ as in the preceding paragraph such that $|f(z)| = |p(z)|$ for
all $z \in {\bf C}$ with $|z| = 1$.

	Let us restrict our attention to complex polynomials as in
(\ref{p(z) = a_n z^n + cdots + a_1 z + a_0}).  If $p(z)$ is as in
(\ref{p(z) = a_n z^n + cdots + a_1 z + a_0}), put
\begin{equation}
	\|p\| = \sup \{|p(z)| : z \in {\bf C}, \ |z| = 1\}
\end{equation}
and
\begin{equation}
	\|p\|_1 = \sum_{j=0}^n |a_j|.
\end{equation}
These define norms on the complex vector space of polynomials, and we
have that
\begin{equation}
	\|p\| \le \|p\|_1.
\end{equation}

	If $p$, $q$ are complex polynomials, then
\begin{equation}
	\|p \, q\| \le \|p\| \, \|q\|
\end{equation}
and
\begin{equation}
	\|p \, q\|_1 \le \|p\|_1 \, \|q\|_1.
\end{equation}
If $p(z)$ is as in (\ref{p(z) = a_n z^n + cdots + a_1 z + a_0}), then
\begin{equation}
	\sum_{j=0}^n |a_j|^2 = \frac{1}{2 \pi} \int_{\bf T} |p(z)|^2 \, |dz|
		\le \|p\|^2,
\end{equation}
where ${\bf T}$ denotes the unit circle in ${\bf C}$.  One can use
this to estimate $\|p\|_1$ in terms of $\|p\|$ and $\|p'\|$,
where
\begin{equation}
	p'(z) = n \, a_n \, z^{n-1} + \cdots + a_1
\end{equation}
is the derivative of $p(z)$.  As a result one can show that
\begin{equation}
	\|p\| = \lim_{l \to \infty} \|p^l\|_1^{1/l}.
\end{equation}

\section{An integral operator}
\label{integral operator}
\setcounter{equation}{0}

	Let $V$ denote the vector space of continuous complex-valued
functions on the unit interval $[0, 1]$ in the real line.  If $f \in
V$, then we put
\begin{equation}
	\|f\| = \sup \{|f(x)| : 0 \le x \le 1\},
\end{equation}
which is the usual supremum norm of $f$.  Define a linear operator $T$
on $V$ by
\begin{equation}
	T(f)(x) = \int_0^x f(s) \, ds.
\end{equation}
If $f$ happens to be real-valued, then $T(f)$ is real-valued, and if
$f(x) \ge 0$ for all $x \in [0, 1]$ too, then $T(f)(x) \ge 0$ for all
$x \in [0, 1]$ as well.  Notice that
\begin{equation}
	\|T(f)\| \le \|f\|
\end{equation}
for all $f \in V$, and that equality holds when $f$ is the constant
function equal to $1$.

	For each positive integer $n$ let $T^n$ denote the $n$-fold
composition of $T$ on $V$.  Equivalently, this is equal to $T$ when $n
= 1$, and in general $T^{n+1}(f) = T(T^n(f))$.  One can express
$T^n(f)$ explicitly as an $n$-fold integral of $f$, and observe that
$T^n(f)$ is real-valued when $f$ is real-valued and nonnegative when
$f$ is nonnegative.  If $f$ is the constant function equal to $1$,
then $T^n(f)(1) = 1/n!$.  Indeed, $T^n(f)(1)$ is equal to the volume
of the points $x = (x_1, \ldots, x_n)$ in ${\bf R}^n$ such that $0 \le
x_1 \le x_2 \le \cdots \le x_n \le 1$.  Using this one can check that
if $f$ is the constant function equal to $1$, then $\|T^n(f)\| =
1/n!$.  For any function $f \in V$ we have that $\|T^n(f)\| \le
\frac{1}{n!} \|f\|$.

	For each $f \in V$ we also have that
\begin{equation}
	\|T(f)\| \le \int_0^1 |f(y)| \, dy.
\end{equation}
As a result, if $f_1, \ldots, f_l$ are elements of $V$, then
\begin{equation}
	\sum_{j=1}^l \|T(f_j)\| \le \biggl\|\sum_{j=1}^l |f_j| \biggr\|.
\end{equation}

\end{document}